\documentclass[a4paper,12pt,leqno]{article}
\usepackage{latexsym}
\usepackage[all]{xy}

\usepackage{amssymb} 
\usepackage{amsmath} 
\usepackage{theorem}

\def\P{{\mathbf{P}}}
\def\Z{{\mathbb{Z}}}
\def\K{{\mathbb{K}}}
\def\R{{\mathbb{R}}}
\def\C{{\mathbb{C}}}
\def\A{{\mathcal{A}}}
\def\B{{\mathcal{B}}}
\def\F{{\mathbb{F}}}

\DeclareMathOperator{\codim}{codim}

\DeclareMathOperator{\Der}{Der}

\numberwithin{equation}{section}

\newcommand{\owari}{\hfill$\square$}

\theoremstyle{break}
\newtheorem{theorem}{Theorem}[section]
\newtheorem{prop}[theorem]{Proposition}
\newtheorem{cor}[theorem]{Corollary}
\newtheorem{lemma}[theorem]{Lemma}

\newtheorem{rem}[theorem]{Remark}

\newtheorem{problem}[theorem]{Problem}
\newtheorem{conj}[theorem]{Conjecture}

\title{Double points of free projective line arrangements}
\author{
Takuro Abe
\footnote
{
Institute of Mathematics for Industry,
Kyushu University,
Fukuoka 819-0395, Japan.
Email:abe@imi.kyushu-u.ac.jp
}
}
\date{\today} 

\begin{document}

\maketitle

\begin{abstract}
We prove Anzis and Toh\u aneanu conjecture, that is the Dirac-Motzkin conjecture for supersolvable line arrangements in the projective plane 
over an arbitrary field of characteristic zero. Moreover, 
we show that a divisionally free arrangements of lines contain at least one double point, that can be regarded as the Sylvester-Gallai theorem for some free arrangements. 
This is a corollary of a general result that if you add a line to a free projective line arrangement, then that line has to contain at least one double point. 
Also we prove some conjectures and one open problems related to 
supersolvable line arrangements and the number of double points. 
\end{abstract}

\section{Introduction}

Let $\K$ be a field, $V=\K^3,\ 
S=\mbox{Sym}^*(V^*)=\K[x,y,z],\ 
\P^2:=\mbox{Proj}(S)$ and 
$\Der S:=S \partial_x\oplus
S \partial_y\oplus
S \partial_z$. Then $\Der S$ is an $S$-graded module and we say that $\theta \in \Der S$ is \textbf{homogeneos} of degree $d$ if $\theta(x),\theta(y),\theta(z)$ belong to the homogeneous part $S_d$ of $S$ of degree $d$.
Let $\A$ be a \textbf{line arrangement} in $\P^2$ (called a \textbf{projective line arrangement}), equivalently, an arrangement 
of linear planes in $V$. For each $H \in \A$ let 
$\alpha_H$ be the defining linear form. 
Let $L(\A)$ be the \textbf{intersection lattice} of $\A$ defined by 
$$
L(\A):=\{ \cap_{H \in \B} H \mid \B \subset \A\}
$$
where we consider everything in $V=\cap_{H \in \emptyset}H$. So $V$ is always contained in $L(\A)$. Also \textbf{we always assume that the origin is in $L(\A)$} (equivalently, $\A$ is \textbf{essential}). Let $L_2(\A)$ be the set $X\in L(\A)$ 
such that $\codim_V X=2$. So they are points in $\P^2$, and we call 
$q \in L_2(\A)$ a \textbf{point}. Define the \textbf{M\"obius function} $\mu:L_2(\A) \rightarrow \Z$ 
by  
$\mu(p):=|\{H \in 
\A \mid p \in H\}|-1$ for $p \in L_2(\A)$. We say that $p \in L_2(\A)$ is a \textbf{double point} if 
$\mu(p)=1$. Let 
\begin{eqnarray*}
n_2(\A):&=&|\{ p \in L_2(\A) \mid \mu(p)=1\}|,\\
n_2(H):&=&|\{ p \in L_2(\A) \mid \mu(p)=1,\ p \in H\}|\ (H \in \A).
\end{eqnarray*}
Then we can define a \textbf{characteristic polynomial} $\chi_0(\A;t)$ of a projective line 
arrangement $\A$ by 
$$
\chi_0(\A;t):=t^2-(|\A|-1)t+
(\sum_{q \in L_2(\A)} \mu(q))-|\A|+1.
$$
Now we can define the \textbf{logarithmic vector field} $D(\A)$ of $\A$ as follows: 
$$
D(\A)=\{\theta \in \Der S \mid \theta(\alpha_H) \in S\alpha_H\ (\forall H \in \A)\}.
$$
$D(\A)$ is a reflexive $S$-graded module of rank $3$, and not free in general. 
We say that $\A$ is \textbf{free} with \textbf{exponents} 
$\exp(\A)=(1,d_2,d_3)$ if there are $\theta_2,\theta_3 \in D(\A)$ of degrees $d_2,d_3$ such that $\theta_2,\theta_3$ and the 
Euler derivation $\theta_E=x\partial_x+wy \partial_y+z\partial_z$ form a free basis for $D(\A)$. By Terao's factorization theorem in \cite{Te}, in this case
$$
\chi_0(\A;t)=(t-d_2)(t-d_3).
$$
Hence in this case 
$|\A|=1+d_2+d_3$. Free arrangements have been intensively studied in the theory of hyperplane arrangements. 
%
Moreover, in the free arrangements, a very important class of arrangements, so called the 
\textbf{supersolvable} arrangements exists. We say that $\A$ is supersolvable if 
there is a point $p \in L_2(\A)$ (called the \textbf{modular point}) such that for $\A_p:=\{H \in \A \mid p \in H\}$ and 
for all distinct $H,L \in \A \setminus \A_p$, there is the 
unique $K \in \A_p$ such that $H \cap L \subset K$. If $\mu(p)=m-1$, then $\A$ is free with 
$\exp(\A)=(1,m-1,|\A|-m)$ (see \cite{OT}, Theorem 4.58 for 
example). For general reference on the above results, please refer 
\cite{OT} and \cite{Y}.

Supersolvable line arrangements have been intensively studied from viewpoints of algebra, combinatorics, geometry, and singularity. 
Among them, one of the 
most interesting and important conjetures is the following.

\begin{conj}[\cite{AnT}, Conjecture 3.2]
Let $\K=\C$, and let $\A$ be supersolvable. Then $n_2(\A) \ge \frac{|\A|}{2}$. 
\label{ATconj}
\end{conj}

Let us review a history of 
Conjecture \ref{ATconj}. 
The origin of Conjecture \ref{ATconj} is the \textbf{Sylvester's problem} asking all 
real line arrangements that is not a pencil in 
$\P^2_\R$ have at least one double points in \cite{Sy}. Here we say that a projective line arrangement $\A$ is \textbf{pencil} if 
all the lines in $\A$ intersects at the same one point. This problem is solved by Gallai in \cite{G} in 1944, thus it is now called the \textbf{Sylvester-Gallai theorem}. 
So it is natural to ask a lower bound of the cardinality of double points for such a line arrangement.
When $\A$ is an arbitrary arrangement and $\K=\R$, Conjecture \ref{ATconj} is 
called the \textbf{Dirac-Motzkin conjecture}, and it is proved to be true when $|\A|$ is sufficiently large by Green and Tao in \cite{GT}. In \cite{AnT}, when $\A$ is supersolvable and 
$\K=\R$, this is proved without the assumption on $|\A|$. However, it is known that 
there is a line arrangement in the complex projective plane which has no double points like the dual Hesse arrangement, see 
\cite{D}, \cite{G}, and \cite{K}. So it is natural to ask the supersolvable version of 
the Dirac-Motzkin conjecture over the complex number field as in Conjectre \ref{ATconj}. 
Also, a supersolvable version of the Sylvester-Gallai theorem is reformulated as follows, which is a 
weaker version of Conjecture \ref{ATconj}. 

\begin{conj}[\cite{HaHa}, Conjecture 11]
Let $\K=\C$, and let $\A$ be supersolvable. Then $n_2(\A) >0$. 
\label{HHconj}
\end{conj}

These conjectures motivated several related researches on supersolvable 
line arrangements form several points of view, e.g., algebra, algebraic geometry, topology and combinatorics aiming at the classification of such arrangements. 
A part of such studies are in \cite{AD}, \cite{AnT}, 
\cite{CHMN}, \cite{DMO}, \cite{DStNSS}, \cite{HaHa}, and \cite{T2}.
On Conjectures \ref{ATconj} and  \ref{HHconj}, there have been several progresses and partial answers in \cite{AD}, \cite{AnT}, \cite{HaHa} (e.g., 
see Theorem \ref{AD}). The main result in this article is to prove (a) Conjecture \ref{ATconj} 
over an arbitrary field $\K$ of characteristic zero in full generality, and (b) 
Conjecture \ref{HHconj} in a wider category of supersolvable arrangements, so called the 
divisionally free arrangements. Let us state the first main result in this article.

\begin{theorem}
Let $\K$ be a field of characteristic zero, 
and let $\A$ be a supersolvable line arrangement in $\P^2$. 
Then $n_2(\A) \ge \frac{|\A|}{2}$. 
\label{main}
\end{theorem}

Actually, we can give a new lower 
bound for $n_2(\A)$ when $\A$ is supersolvable as follows.

\begin{cor}
Let $\A$ be a line arrangement in $\P^2$ over a field $\K$ of characteristic zero. 
Assume that $\A$ is supersolvable with $|\A|=m+k,\ k \ge 1$, $p \in L_2(\A)$ a modular point with 
$\mu(p)=m-1 \ge 1$. Note that $\exp(\A)=(1,m-1,k)$. Then we have the following:
\begin{itemize}
\item[(1)]
If $k \le m$, then $n_2(\A) \ge k(m-k+1)$.
\item[(2)]
If $k \ge m$, then $n_2(\A) \ge k$.
\end{itemize}
\label{bound}
\end{cor}

\begin{rem}
If $k \le m$, then a lower bound 
$$
n_2(\A) \ge k(m-k+1) \ge \frac{|\A|}{2}
$$
is proved in Theorem 1.7 in \cite{AD} when $\K=\C$. In this case 
Conjecture \ref{ATconj} holds true. We give another proof of this inequality over an 
arbitrary field $\K$ of characteristic zero in the proof of Theorem \ref{main}.
\label{ADbound}
\end{rem}

Next let us state the second main result. 
For that, let us introduce a \textbf{divisionally free arrangements} of lines. 
We say that $\A$ is divisionally free if 
there is $H \in \A$ such that 
$\chi_0(\A;|\A^H|-1)=0.
$
Here $\A^H:=\{ 
q \in L_2(\A) \mid q \in H\}$. Note that divisional freeness is a combinatorial property, and divisionally free arrangements contain the famous inductively free arrangements. Thus in this class, Terao's conjecture is true, that asserts that the freenees is combinatorial. Supersolvable arrangements are divisionally free, see \cite{A} and \cite{A2} for details of divisional freeness. In this class we can show the Sylvester-Gallai theorem over an arbitrary field of characteristic zero.

\begin{theorem}
If $\A$ is a 
divisionally free line arrangement in $\P_\K^2$ over an arbitrary field of characteristic zero, then $n_2(\A)>0$.
\label{main3}
\end{theorem}

Actually, Theorem \ref{main3} is a corollary of the following general result, that plays an important role in the proof of Theorem \ref{main}.

\begin{theorem}
Let $\A$ be a projective line arrangement over an arbitrary field of characteristic zero $\K$, and $H \in \A$. If $\A':=\A\setminus \{H\}$ is free, then $n_2(\A)>n_2(H)>0$. 
Equivalently, if you add a line to a free projective line arrangement, then that line has to contain at least one double points.
\label{main10}
\end{theorem}

In Theorem \ref{main10}, $\A$ could be not free. Also note that Theorem \ref{main10} fails if $\K$ has 
a positive characteristic. See Remark \ref{chap} for details.
\medskip

The organization of this 
article is as follows. In \S2 we introduce several results and definitions for the proof. 
In \S3 we prove main results posed in 
this section, and also prove several conjectures and open problems by using them. In \S4 we pose several conjectures related to the cardinality of double points of non-free projective line arrangements.
\medskip

\noindent
\textbf{Acknowledgements}. The author is 
supported by JSPS KAKENHI 
Grant-in-Aid for Scientific Research (B) Grant Number 
JP16H03924. The author is grateful to Hiraku 
Kawanoue for his advice and informing Proposition 
\ref{kawanoue} to the author.

\section{Preliminaries}

\textbf{From now on we assume that $\K$ is a field of characteristic zero} unless otherwise 
specified. Our approach to prove them is purely algebraic, and different from the previous approaches. 
To prove Theorem \ref{main} we need a few results and definitions. First let us recall the multiarrangement 
on $\K^2$. Let $\A$ be a central arrangement of lines in $V':=\K^2$ and let 
$\mathbf{m}:\A \rightarrow \Z_{> 0}$ be a multiplicity. Then a pair 
$(\A,\mathbf{m})$ is called a 
\textbf{multiarrangement}. Let $S':=\K[x,y]$ be the coordinate ring of $V'$. Then we can define the logarithmic derivation module $D(\A,\mathbf{m})$ as 
$$
D(\A,\mathbf{m}):=\{
\theta \in \Der S' \mid 
\theta(\alpha_H) \in S'\alpha_H^{\mathbf{m}(H)}\ (\forall 
H \in \A)\}.
$$
It is known (see \cite{Z} for example) that 
$D(\A,\mathbf{m})$ is always free, so we have its exponents $\exp(\A,\mathbf{m})=(d_1,d_2)$. It is known that 
$d_1+d_2=|\mathbf{m}|:=\sum_{H \in \A} \mathbf{m}(H)$. 
By definition, we have the following easy but 
important lemma.

\begin{lemma}
Let $\A_1 \subset \A_2$ be arrangements in $\K^2$ and 
$\mathbf{m}_i:\A_i \rightarrow \Z$ be multiplicities. Assume that 
$\mathbf{m}_2(H) \ge \mathbf{m}_1(H)$ for all $ 
H \in \A_1$. Then $D(\A_2,\mathbf{m}_2) \subset D(\A_1,\mathbf{m}_1)$.  

In particular, for $\exp(\A_1,\mathbf{m}_1)=(d_1,d_2)$ and $\exp(\A_2,\mathbf{m}_2)=(e_1,e_2)$ with $d_1 \le d_2,\ e_1 \le e_2$, it holds that 
$d_1 \le e_1$ and $d_2 \le e_2$.
\label{key}
\end{lemma}

It is not easy to determine exponents of multiarrangements, but in some cases we 
can determine them completely. The following is one of them and it plays the 
key role in our proof.

\begin{prop}[\cite{ST}, Proposition 5.4. See \cite{WY}, Example 2.2 too]
Let $\A$ be a line arrangement in $\K^2$ and 
$\mathbf{2}$ be the constant multiplicity $2$ on $\A$, i.e., 
$\mathbf{2}(H)=2$ 
for any $H \in \A$.
Then $\exp(\A,\mathbf{2})=(|\A|,|\A|)$.
\label{WY}
\end{prop}

\begin{rem}
Note that 
Proposition \ref{WY} holds true only over the field of characteristic zero. For example, consider 
the multiarrangement $(\A,\mathbf{m})$ in $\F_2^2$ defined by $x^2y^2(x-y)^2=0$. 
Then $D(\A,\mathbf{m})$ has a basis 
$x^2 \partial_x+y^2\partial_y,
x^4 \partial_x+y^4\partial_y$. Thus $\exp(\A,\mathbf{m})=(2,4)$.
\label{rem1}
\end{rem}

From a line arrangement $\A$ in $\P^2$ and $H \in \A$, 
we can canonically construct a multiarrangement 
$(\A^H,\mathbf{m}^H)$ in $H \simeq \K^2$ as follows: $\A^H:=\{H \cap L 
\mid L \in \A \setminus \{H\}\}$, and 
$$
\mathbf{m}^H(X):=|\{
L \in \A \setminus \{H\} \mid 
L \cap H=X\}|
$$
for $X \in \A^H$. Here $\A^H$ is an arrangement in $H \simeq \K^2$. 
To relate $D(\A)$ and $D(\A^H,\mathbf{m}^H)$, we need an $S$-graded submodule $D_H(\A) \subset D(\A)$. Namely, 
for 
$H \in \A$, the $S$-graded submodule $D_H(\A)$ is defined by 
$$
D_H(\A):=\{\theta \in D(\A) \mid \theta(\alpha_H)=0\}.
$$
It is known that $D(\A)=S\theta_E \oplus D_H(\A)$ for any $H \in \A$ (e.g., see Lemma 1.33, \cite{Y}). 
So $D_H(\A) \simeq D_L(\A)$ for any lines $H, L$ in $\A$. Then by \cite{Z}, there is the 
\textbf{Ziegler restriction map} $\pi_H:D_H(\A) 
\rightarrow D(\A^H,\mathbf{m}^H)$ defined by 
$$
\pi_H(\theta)(\overline{f}):=\overline{\theta(f)},
$$
where for 
$f \in S$, $\overline{f}$ indicates the image of $f$ in $S/\alpha_H S$. 
Then Ziegler proved the following.

\begin{theorem}[\cite{Z}]
Let $H \in \A$. Then 

(1)\,\,
there is an exact sequence 
$$
0 \rightarrow \alpha_H D_H(\A) \rightarrow D_H(\A) 
\stackrel{\pi_H}{\rightarrow} 
D(\A^H,\mathbf{m}^H).
$$
(2)\,\,
Assume that $\A$ is free with $\exp(\A)=(1,d_2,d_3)$ and let $H \in \A$. Then $\exp(\A^H,\mathbf{m}^H)=(d_2,d_3)$.
\label{Z}
\end{theorem}

Next let us recall some fundamental results on supersolvable line arrangements. 

\begin{lemma}
Let $\A$ 
be a supersolvable line arrangement with 
a modular point $p \in L_2(\A)$ such that $\mu(\A)=m-1$ and $|\A|=m+k$. Thus 
$\exp(\A)=(1,m-1,k)$. 

(1)\,\,
Let $H \in \A \setminus \A_p$. Then $
|\A^H|=m$.

(2)\,\,
Let $\A \setminus \A_p=\{H_1,\ldots,H_k\}$. Then 
$
n_2(\A) \ge \sum_{i=1}^k n_2(H_i)$.
\label{1}
\end{lemma}

\noindent
\textbf{Proof}.
(1)\,\, Since $|\A_p|=m$ and $H \not \in \A_p$, it is clear that $|\A^H|\ge m$.
Let $q \in L_2(\A)$ be a point on $H$ that is not on a line in $\A_p$. Then 
$q=H \cap L$ for $L \in \A \setminus (\A_p \cup\{H\})$. Since 
$\A$ is supersolvable and $p$ is a modular point, 
there is $K \in \A_p$ such that $H \cap L =q \in K$, a contradiction. So  $|\A^H|=m$.

(2)\,\, 
Let $q_1 \in H_1$ and $q_2 \in H_2$ be double points 
belonging to $L_2(\A)$. It suffices to show that $q_1 \neq q_2$. Assume that $
q_1=q_2$. Then $q_1=q_2=H_1 \cap H_2$. 
Since 
$\A$ is supersolvable and $p$ is a modular point, 
there is $K \in \A_p$ such that $H_1 \cap H_2 =q_1=q_2 \in K$. So 
$q_1=q_2=H_1\cap H_2\cap K$. As a 
consequence, $q_1=q_2$ is not a double point, a contradiction. 
\owari
\medskip

To count $|\A^H|$, the following is the key.

\begin{theorem}[\cite{A}, Theorem 1.1, \cite{A2}, Theorem 1.1, Corollary 1.2]
Let a line arrangement $\A$ be free with $\exp(\A)=(1,d_2,d_3)$ with $d_2\le d_3$. Then 
\begin{itemize}
\item[(1)]
$|\A^H|\le d_2+1$ or $|\A^H|=d_3+1$. 
\item[(2)]
For $L \not \in \A$, let $\B:=\A \cup \{L\}$. Then 
$|\B^L|=1+d_2$ or $|\B^L| \ge d_2+d_3+1$.
\item[(3)]
$\A$ is (divisionally) free if $|\A^H| = d_2+1$ or $|\A^H|=d_3+1$.
\end{itemize}
\label{A1}
\end{theorem}

Finally recall a partial result on Conjecture \ref{ATconj}. 

\begin{theorem}[\cite{AD}, Theorem 1.7]
Let $\A$ be a line arrangement in $\P^2$ over $\C$.
Assume that $\A$ is supersolvable with $|\A|=m+k,\ k \ge 1$, $p \in L_2(\A)$ a modular point with 
$\mu(p)=m-1 \ge 1$. 
If $k \le m$, then $n_2(\A) \ge k(m-k+1)$.
\label{AD}
\end{theorem}

\section{Proof of main results}

First we prove Theorem \ref{main10}.
\medskip

\noindent
\textbf{Proof of Theorem \ref{main10}}. 
Let $\exp(\A')=(1,a,b)$ with $a \le b$. 
By Theorem \ref{A1} (2), 
$|\A^H|=1+a$ or at least $b+1$. First assume the latter. 
Since $|\mathbf{m}^H|=a+b+1 
\le 2(b+1)=2|\A^H|$, it holds that $n_2(H)>0$. Thus we may assume that 
$|\A^H|=a+1$, thus $\A$ is free with $\exp(\A)=(1,a,b+1)$ by Theorem \ref{A1} (3). Assume that 
$n_2(H)=0$. Then by definition, $m^H(X) \ge 2$ for all $X \in \A^H$. So Lemma \ref{key} shows that 
$D(\A^H,\mathbf{m}^H) \subset 
D(\A^H,\mathbf{2})$. By Theorem \ref{Z} (2), 
$\exp(\A^H,\mathbf{m}^H)=(a,b+1)$ and by Proposition \ref{WY}, $\exp(\A^H,\mathbf{2})=(
a+1,a+1)$, which contradicts Lemma \ref{key}.\owari
\medskip

\noindent
\textbf{Proof of Theorem \ref{main3}}. 
Combining Theorems \ref{A1} and  
\ref{main10}, we have $n_2(\A)>0$. \owari
\medskip

\begin{rem}
Theorem \ref{main10} 
is not true if $\mbox{ch}(\K) \neq 0
$. Let $\K=\F_2$ and let $\A$ be a line arrangement in $\P^2_{\F_2}$ consisting of all 
lines in $\P^2_{\F_2}$. Thus $|\A|=7$. Let $\A \ni H:x=0$. Then it is easy to show that 
$D(\A \setminus \{H\})$ is free with basis 
$$
\theta_E,\ x^2\partial_x+y^2\partial_y+z^2\partial_z,\ 
(x+y)(x+z)(x+y+z)\partial_x.
$$
However, it is also clear that $n_2(H)=0$. 

Of course even if $\mbox{ch}(\K)=p>0$, it is easy to show that Theorem \ref{main10} is true 
if (1) $\A$ is not free in terms of Theorem \ref{main10}, or (2) 
$a$ is not divisible by $p$. See the proof of Proposition 1.53 in \cite{Y}.
\label{chap}
\end{rem}

Now we can prove 
Theorem \ref{main}. 
\medskip

\noindent
\textbf{Proof of Theorem \ref{main}}. 
Let $p \in L_2(\A)$ be the modular point with $\mu(p)=m-1$ and let $|\A|=k+m$. Thus 
$\exp(\A)=(1,m-1,k)$. If $m=1$, then the statement is clear. So we may assume that 
$m \ge 2$. First we prove when $ k\le m$. If $k=0$, then 
it is not essential. So we may assume that $k \ge 1$. 
If $k=1$, then it is the coning of $\A_p$, so $n_2=m$, which satisfies the inequality. 
Let $k \ge 2$, and let $H \in \A \setminus \A_p$, and count the number of 
double points $n_2(H)$ on $H$. 
By Lemma \ref{1} (1), $|\A^H|=m$. 
Let $q \in L_2(\A) \cap \A^H$ be a non-double point. Since $p$ is a modular point and 
$H \not \in \A_p$, there is $L \in \A \setminus (\A_p \cup\{H\})$ such that 
$H \cap L=q$. Since 
$\A$ is supersolvable again, there is $K \in \A_p$ such that $K \ni q=H \cap L \cap K$. Thus $|\A^H|-n_2(H)
\le |\A \setminus (\A_p \cup \{H\})|=k-1$. 
Thus $n_2(H) \ge m-(k-1)>0$ since $m \ge k$. Now Lemma \ref{1} (2) implies that 
$n_2(\A) \ge (m-k+1)k$. Now compare
\begin{eqnarray*}
2(n_2(\A)-\frac{|\A|}{2})&\ge&
2k(m-k+1)-(m+k)\\
&=&(m-k)(2k-1) \ge 0.
\end{eqnarray*}
So the proof is completed when $ k \le m$.
%
%
%

Next consider when $k \ge m$. The statement is clear if $|\A|-|\A_p|=1$. So 
$|\A|-|\A_p|\ge 2$. 
Let $H \in \A \setminus \A_p$. By definition of the supersolvable arrangement, 
$\A \setminus \{H\}$ is supersolvable, thus free. 
Hence $n_2(H)>0$ by Theorem \ref{main10}.
%
Now let $\A \setminus \A_p=\{H_1,\ldots,H_k\}$. 
Since $n_2(H_i)>0$ for $i=1,
\ldots,k$, Lemma \ref{1} (2) shows 
that 
$$
n_2(\A) \ge \sum_{i=1}^k n_2(H_i) \ge 
k.
$$
Now compute 
\begin{eqnarray*}
2(n_2(\A)-\frac{|\A|}{2})&\ge&
2k - (m+k)=k-m \ge0, 
\end{eqnarray*}
which completes the proof. \owari
\medskip

\noindent
\textbf{Proof of Corollary \ref{bound}}. 
Immediate by the proof of Theorem 
\ref{main}. \owari
\medskip


Based on the proof of Theorem \ref{main}, we can show the following
general relation between $n_2(H)$ with $H \in \A$ and 
the \textbf{minimal degree relation} $r(\A)=mdr(\A)$ defined by
$$
r(\A)=mdr(\A):=\min\{d \mid D_H(\A)_d \neq (0)\}
$$
here $M_d$ for an $S$-graded module $M$ indicates the homogeneous degree $d$-part of $M$. 
Then the most general version of Theorem \ref{main10} is given as follows:

\begin{theorem}
Let $\A$ be a projective line arrangement with $H \in \A$ such that 
$r(\A)<|\A^H|$. Then $n_2(\A)>n_2(H)>0$. 
\label{HaHa2}
\end{theorem}

\noindent
\textbf{Proof}. Apply Theorem \ref{WY} and the same argument in the proof of Theorem 
\ref{main10}. \owari
\medskip

Thus in the terminology of $mdr(\A)=r(\A)$, we can give a lower bound of $n_2(\A)$. 

\begin{theorem}
Let $\A$ be a line arrangement in $\P^2$ with $r=mdr(\A)$. 
If there is $H \in \A$ with $|\A^H|=k>r$, then $n_2(H) \ge k-r$. 
Moreover, if $\A_{> r}:=\{L \in \A \mid 
|\A^L|> r\}$, then 
$$
n_2(\A) \ge \displaystyle \frac{1}{2} \sum_{L \in \A_{>r}}
(|\A^L|-r).
$$
In particular, 
$$
n_2(\A) \ge \displaystyle \frac{|\A_{>r}|}{2},
$$ 
and $n_2(\A)>0$ if $\A_{>r} \neq 
\emptyset$.
\label{main2}
\end{theorem}

\noindent
\textbf{Proof}. 
Assume that $n_2(H) <k-r$. Let 
$(\B,\mathbf{m})$ be a  multiarrangement of lines in $\K^2$ obtained by removing all 
double points on $\A^H$, and 
putting $\mathbf{m}:=\mathbf{m}^H|_{\B}$. Thus $|\B| = k-n_2(H)>r$, and 
$\mathbf{m}(X) \ge 2$ for all $X \in \B$. 
By definition of $r$, there is $0 \neq \theta \in D_H(\A)_r$. Assume that $\theta$ is divisible by $\alpha_H$. Then by definition of $D_H(\A)$, it holds that $\theta/\alpha_H \in D_H(\A)_{r-1}=(0)$, a contradiction. 
Thus $\theta \not \in \alpha_H 
D_H(\A)_{r}$. Since $\ker (\pi_H)=\alpha_H D_H(\A)$ by Theorem \ref{Z} (1), it holds that 
$0 \neq 
\pi_H(\theta) \in D(\A^H,\mathbf{m}^H)_{r} \subset D(\B,\mathbf{m})_{r} 
\subset D(\B,\mathbf{2})_{r}=(0)$ because $\exp(\B,\mathbf{2})=(k-n_2(H),k-n_2(H))$ with 
$k-n_2(H)>r$ by 
Lemma \ref{key} and Proposition \ref{WY}, a contradiction. The rest 
statements are clear. \owari
\medskip

\begin{rem}
Theorem \ref{main2} is weaker than the argument used in the proof of Theorem \ref{main} 
since all double point in the proof of Theorem \ref{main} are distinct, but those in Theorem \ref{main2} are not. The difference comes from whether there is a modular point or not.
\end{rem}



In \cite{AnT}, two equivalent Conjectures 3.3 and 3.6 are posed over $\C$, which are stronger than Conjecture \ref{ATconj}. 
By Theorem \ref{main2} or the same proof as in Theorem \ref{main}, we can prove them (explicitly, we can show \cite{AnT}, Conjectures 3.3 as Theorem \ref{AT2}, and Conjecture 3.6 follows by the equivalence) over an arbitrary field of charactersitic zero.

\begin{theorem}[\cite{AnT}, Conjectures 3.3]
There are no supersolvable line arrangement $\A$ in $\P^2$ with a modular point $p$ of 
$\mu(p) \ge 2$ such that 
$n_2(H)=0$ for some $H \in \A \setminus \A_p$.
\label{AT2}
\end{theorem}

\begin{theorem}[\cite{AnT}, Conjectures 3.6]
Let $\K$ be an arbitrary field of characteristic zero, $f_1,\ldots,f_n \in \K[x,y]$ be $m$-linear forms with $\mbox{gcd}(f_i,f_j)=1$ for 
any $i \neq j$. Then 
$$
\prod_{i=1}^n f_i \not \in \bigcap_{1 \le i < j \le m} \langle f_i+z,f_j+z\rangle,
$$
where 
$\langle f_i+z,f_j+z\rangle$ is the ideal in $\K[x,y,z]$ generated by $f_i+z$ and 
$f_j+z$.
\label{AT3}
\end{theorem}

Moreover, Question 16 posed in \cite{HaHa} can be settled affirmatively as follows:

\begin{theorem}
Let $\A$ be a supersolvable line arrangement in $\P^2$ over a filed of characteristic zero. Let $p$ be a modular point of $\A$
such that 
$\mu(p) \ge 1$.
Assume that $|\A \setminus \A_p| \ge 1$. Then 
$n_2(\A) \ge \max\{|\A|-m,m\}$.
\label{large}
\end{theorem}

\noindent
\textbf{Proof}. Let $|\A|=m+k$, so $\exp(\A)=(1,m-1,k)$. By the assumption, 
$k \ge 2$. 
First assume that $m \ge k$. Then Corollary \ref{bound} shows that 
$n_2(\A) \ge k(m-k+1)$, and 
$$
k(m-k+1)-m=k (m-k)-(m-k) \ge 0
$$
since $k \ge 1$. Thus 
$n_2(\A) \ge m=\max\{m,|\A|-m=k\}$. Assume that $m\le k$. Then $n_2(\A) \ge k$ by 
Corollary \ref{bound}. Thus 
$n_2(\A) \ge k=\max\{m,|\A|-m=k\}$. \owari
\medskip

As mentioned in \cite{HaHa}, Theorem \ref{large} implies 
Conjecture \ref{ATconj} too. 
\medskip

\section{Conjectures}

Based on the previous results, the following problem is natural to ask.

\begin{problem}
Consider the Sylvester's problem depending on $D(\A)$. For example, consider it when $\A$ is free.
\label{Sylvester}
\end{problem}

For a generic line arrangement $\A$ in $\P_\K^2$, the Dirac-Motzkin conjecture is true. Also, when 
we construct a free line arrangement $\A$, experimentally, we know that $n_2(\A)$ decreases. 
Based on these very rough observation with analysis of the known Sylvester-Gallai configurations, we pose the following conjecture.

\begin{conj}
Let $\A$ be a projective line arrangement in $\P^2_\K$, where $\K$ is a field of characteristic zero. 
If $\A$ is not free, then
$n_2(\A)>0$.
\label{conj1}
\end{conj}

A weaker version is as follows:

\begin{conj}
Let $\A$ be a free projective line arrangement in $\P^2_\K$, where $\K$ is a field of characteristic zero. 
Let $H \in \A$. Assume that $\A\setminus \{H\}$ is not free. Then
$n_2(\A\setminus \{H\})>0$.
\label{conj2}
\end{conj}

Conjecture \ref{conj2} is a deletion version of Theorem \ref{main10}, i.e., Theorem \ref{main10} asserts that 
any addition to a free arrangement has a double point. Then how about the deletion? Since the dual Hesse arrangement, that consists of $9$-lines with $12$-triple points and no other intersection points, can be 
obtained by deleting a line from a free arrangement (see \cite{ACKN}, Theorem 1.1), Conjecture \ref{conj2} is not true if $\A \setminus \{H\}$ is free. Even in that case, $D(\A \setminus \{H\})$ has a good algebraic structure so called the \textbf{plus-one generated} property (see \cite{A3}). Thus we may have a chance to approach Conjecture \ref{conj2} by using algebraic technique. 
Actually, the supersolvable version of Conjecture \ref{conj2} is true as follows:

\begin{theorem}
Let $\A$ be a supersolvable projective line arrangement with a modular point $p \in L_2(\A)$, $
\mu(p)=m-1\ge 1$. Let $\A':=\A \setminus \{H\}$. Assume that $\A'$ is not a pencil. 
Then $n_2(\A') >0$.
\label{ns}
\end{theorem}

\noindent
\textbf{Proof}. 
Assume that $H \in \A \setminus \A_p$. Since $\A'$ is not a pencil, $\A'$ is still supersolvable. Thus Theorem \ref{main} shows that $n_2(\A') \ge |\A'|/2>0$. Next assume that 
$H \in \A_p$. By Theorem \ref{main10}, there is at least one double point on any $L \in \A 
\setminus \A_p$. Assume that $n_2(\A')=0$, that occurs only when all such double points are on $H$. 
So $|\A^H|=1+|\A|-m$. Hence it is clear that $\A'$ is still supersolvable with a modular point $p$, $\mu(p)=m-2$, contradicting Theorem \ref{main}. Thus $n_2(\A')>0$. 
\owari
\medskip

\begin{rem}
Note that such an $\A'$ in Theorem \ref{ns} is called \textbf{nearly supersolvable} in 
\cite{DStNSS} if $\A'$ is not supersolvable. Thus the proof of Theorem \ref{ns} asserts that 
if $\A$ is nearly supersolvable, then $n_2(\A)>0$.
\end{rem}

So far Conjecture \ref{conj1} is true in the following class:

\begin{theorem}
For a projective line arrangement $\A$, $n_2(\A)>0$ if
\begin{itemize}
\item[(1)]
$\A$ has a generic line (thus not free),
\item[(2)]
$\A$ is (nearly) supersolvable,
\item[(3)]
there is $H \in \A$ such that $\A \setminus \{H\}$ is free (in particular, inductively and divisionally free arrangements),
\item[(4)]
there is $H \in \A$ and $p \in \A^H$ such that 
$|m^H(p)| \ge |\A|/2$, or 
\item[(5)]
there is $H \in \A$ such that $|\A^H| >mdr(\A)=r(\A)$.
\end{itemize}
\end{theorem}

\noindent
\textbf{Proof}. Immediate by Theorems \ref{main}, \ref{main10}, \ref{HaHa2}, \ref{ns}, and the proof of Theorem \ref{main10}.\owari
\medskip

It is natural to ask whether the Dirac-Motzkin conjecture holds for non-free arrangements. 
However, it is pointed it out by Hiraku Kawanoue that this is not true even if $\A$ 
can be obtained by deleting a line from a free arrangement in general.

\begin{prop}[Kawanoue]
There is a free line arrangement $\A$ in $\P^2_\C$ such that $\A \setminus \{H\}$ is not free, and 
$n_2(\A \setminus \{H\}) <\frac{|\A|-1}{2}$
for any $H \in \A$.
\label{kawanoue}
\end{prop}

\noindent
\textbf{Proof}. 
Let $\A$ be 
defined by
$$
(x^4-y^4)(y^4-z^4)(x^4-z^4)=0
$$
in $\P_\C^2$. This is called the monomial arrangement with respect to the group 
$G(4,4,3)$, see \cite{ACKN} Theorem 1.1 and 
section 5. It is known that (e.g., see \cite{ACKN}) $\A$ is free with $\exp(A)=(1,5,6)$, 
$|\A^H|=5$, $\A^H$ consists of 
one quadruple points and four triple points 
for any $H \in \A$.
It is also shown in \cite{ACKN} that $\A':=\A \setminus \{H\}$ is not free. Also, it is clear that 
$n_2(\A')=4 <|\A'|/2=11/2$.\owari
\medskip

Though the Dirac-Motzkin conjecture is not true for non-free arrangements in $\P^2_\C$, as in Proposition \ref{kawanoue}, 
Theorem \ref{main} shows that it is true for supersolvable arrangements. So 
let us pose the following problem.

\begin{problem}
Let $\K$ be a field of characteristic zero. Then 
in which class of line arrangements in $\P^2_\K$ the inequality 
$n_2(\A) \ge |\A|/2$ holds?
\label{DiracMotzkin}
\end{problem}

\end{document}